\documentclass[12pt]{article}

\usepackage[T2A]{fontenc}
\usepackage[utf8]{inputenc}

\usepackage{amsfonts}
\usepackage{amssymb}
\usepackage{amsmath}
\usepackage{amsthm}

\def \le {\leqslant}
\def \ge {\geqslant}

\begin{document}

\begin{Large}
\centerline{\bf \"Uber  die rationalen Punkte auf der Sph\"are}
 
\end{Large}
\vskip+1cm
\begin{large}
\centerline{von Nikolay Moshchevitin\footnote{ Die Untersuchung ist  von der Beihilfe der russischen Regierung 
11.
G34.31.0053 und  RFBR  No.12-01-00681-a unterstützt.
 } (Moskau)}
 \end{large}
\vskip+2cm

Wir besch\"aftigen uns hier mit der Approximation von  Punkten auf der  $n$-dimensionalen  Sph\"are durch
 rationale Punkte  der
$n$-dimensionalen  Sph\"are.
Wir geben   einen kurzen Beweis des h\"ubschen Satzes von Kleinbock und Merrill [6], im einfachen Fall $n=2$.

\vskip+0.3cm
{\bf \S 1. 
Parametrisierung der Sph\"are.}

Es sei
$$
\hbox{\got S}^n = \{ {\bf x} = (x_1,...,x_{n+1}) \in \mathbb{R}^{n+1}:\,\, x_1^2+...+x_{n+1}^2 = 1\}
$$
eine Sph\"are  vom Radius $1$ im Euklidischen  Raum $ \mathbb{R}^{n+1}$.
F\"ur 
${\bf t} = (t_1,...,t_n)\in \mathbb{R}^n$
setzen wir die Funktionen
$$
f_j({\bf t} )=
\frac{2t_j}{1+t_1^2+...+t_n^2},\,\, j = 1,...,n,\,\,\,\,\
f_{n+1}({\bf t})=\frac{1-t_1^2-...-t^2_n}{1+t_1^2+...+t_n^2}.
$$
 {Bekanntlich ist,  jeder rationale Punkt 
$$
{\bf x} =  \left(\frac{A_1}{Q},...,\frac{A_{n+1}}{Q}\right),\,\, A_1,...,a_{n+1}, Q \in \mathbb{Z},\,\, (Q,A_1,...,A_{n+1})=1
$$
auf der Sph\"are $\hbox{\got S}^n$ ist von der Form
\begin{equation}\label{ret}
\frac{A_j}{Q} = f_j({\bf t}) =
\frac{2b_jq}{q^2+b_1^2+...+b_n^2},\,\, j = 1,....,n,\,\,\,\,\
\frac{A_{n+1}}{Q} = f_{n+1}({\bf t})=
\frac{q^2-b_1^2-...-b_n^2}{q^2+b_1^2+...+b_n^2}
,
\end{equation}
wo
$$
{\bf t} = \left(
\frac{b_1}{q},...,\frac{b_n}{q}
\right),\,\,\,\,
 q\in \mathbb{Z}_+,b_1,...,b_n \in \mathbb{Z}
,
\,\,\,\,
(q,b_1,...,b_n) = 1.
$$
}
F\"ur
$\alpha = (\alpha_1,...,\alpha_{n+1}) \in \hbox{\got S}^n$ setzen wir
$
\alpha \mapsto \beta = \beta (\alpha ) = (\beta_1,...,\beta_n)$
$\alpha_j = f_j(\beta), \, j = 1,...,n+1$.

\vskip+0.3cm
{\bf \S 2. Resultate.}

E. Hlawka [5] hat folgendes gezeigt.
{Sei $\alpha \in (0,1)$,
dann gibt es zu jedem gen\"ugend  gro\ss en $N>1$, nicht negative Zahlen
$u, v$,
so dass $1\le v \le N$ und
$$
\left|\alpha-\frac{v^2-u^2}{v^2+u^2}\right|\le \frac{2}{Nv},\,\,\,\,
\left|\sqrt{1-\alpha^2}-\frac{2uv}{v^2+u^2}\right|\le \frac{2}{Nv}
$$
gilt.
}

L. Fukshansky  [4] bemerkt, dass  f\"ur $(\alpha_1,\alpha_2) \in \hbox{\got S}^1, \{\alpha_1, \alpha_1\} \not\subset\{0,\pm 1\}$ die Behauptung von Hlawka 
 unendlich viele rationalen Punkten $\left(\frac{A_1}{Q},\frac{A_2}{Q}\right)\in \hbox{\got S}^1$ liefert mit
$$
\max_{i=1,2} \left|\alpha_i - \frac{A_i}{Q}\right|\le \frac{2\sqrt{2}}{Q}.
$$

Kleinbock und Merrill [6] hatten  folgendes gezeigt. 
F\"ur jedes $n$ es gibt eine sehr gro\ss e positive Konstante $C_n$ mit
$$
\min_{
\left(\frac{A_1}{Q},,...,\frac{A_n}{Q}\right)\in \hbox{\got S}^{n+1},\,\,
1\le Q\le T
}\,\,\,\,\,
 \left|\alpha_i - \frac{A_i}{Q}\right| \le \frac{C_n}{\sqrt{QT}}\,\,\,\,\,
\text{f\"ur alle}\,\,\,\,\, T \ge 1,
$$
sodass f\"ur jedes $\alpha \in\hbox{\got S}^{n+1}$  unendlich viele rationale Punkte
$\left(\frac{A_1}{Q},,...,\frac{A_n}{Q}\right)\in \hbox{\got S}^{n+1}$  existieren mit
$$
\max_{1\le i \le n} \left|\alpha_i - \frac{A_i}{Q}\right|\le \frac{C_n}{Q}.
$$

Im Falle $n=1$ wir wollen die folgende Behauptung beweisen:
\vskip+0.3cm
{\bf Satz 1.} \, {\it
Es sei $\alpha =(\alpha_1,\alpha_2) \in \hbox{\got S}^1 \setminus \mathbb{Q}^2$ und $\varepsilon >0$. Es gibt unendlich viele
rationale Vektoren $\left(\frac{A_1}{Q},\frac{A_2}{Q}\right)\in \hbox{\got S}^1\cap \mathbb{Q}^2$, so dass 
\begin{equation}\label{11}
\sqrt{\sum_{i=1,2}
\left(
\alpha_i -\frac{A_i}{Q}
\right)^2}
<\frac{1+\varepsilon}{\sqrt{2}Q}.
\end{equation}

\vskip+0.3cm
{\bf Bemerkung 1.}\,\,{\it
Im Satz 1 kann $\sqrt{2}$ durch keine gr\"o\ss ere Zahl ersetzt werden.}
}
\vskip+0.3cm

 Wir wollen jetzt im Falle $n=2$  die folgenden Behauptungen beweisen:

\vskip+0.3cm
{\bf Satz 2.}
\,
{\it Es sei $\beta_1, \beta_2 \in \mathbb{R}$,
$ T\ge 1$. Es gibt einen  rationalen  Vektor $\left(\frac{b_1}{q}, \frac{b_2}{q}\right)$ 
mit den folgenden Eigenschaften:

\noindent
{\rm (i)} $ 1\le q\le T $;

\noindent
{\rm (ii)}  $\sum_{i=1,2} \left( q\beta_i -{b_i}\right)^2 <\frac{4q}{T}$

\noindent
{\rm (iii)}  $ b_1^2 + b_2^2 \equiv 0 \pmod{q}$.
 
\vskip+0.3cm
{\bf Satz 3.}
\, Sei $\gamma > \sqrt{\frac{3}{\pi}}$.
 Sei $\beta_1$ oder $\beta_2$  nicht rational.
Dann existieren  unendlich viele   rationalen  Vektoren $\left(\frac{b_1}{q}, \frac{b_2}{q}\right)$ 
mit {\rm  (iii)} und

\noindent
{\rm (ii$^*$)}  $\sum_{i=1}^2 \left( q\beta_i -{b_i}\right)^2 <{\gamma}^2$.
 
 }
\vskip+0.3cm
{\bf Bemerkung 2.} 
  Man kann  die Ungleichung (ii)   in  folgende Form
umschreiben:
$$
\sqrt{\sum_{i=1,2}
\left(
\beta_i -\frac{b_i}{q}
\right)^2}
<\frac{2}{\sqrt{qT}}
.$$
 
Aus den  S\"atzen  2,3  und (\ref{ret}) erhalten wir die Folgerungen:

\vskip+0.3cm
{\bf Satz 4.}
\,
{\it Es sei  $\alpha = (\alpha_1,\alpha_2,\alpha_3) \in
\hbox{\got S}^2 $ und
$ T\ge 1$. Es gibt einen  rationale  Vektor $\left(\frac{A_1}{Q}, \frac{A_2}{Q},\frac{A_3}{Q}\right)$ 
mit den folgenden Eigenschaften:

\noindent
{\rm (i)} $ 1\le Q\le T $;

\noindent
{\rm(ii)} $
\sqrt{\sum_{i=1,2,3}
\left(
\alpha_i -\frac{A_i}{Q}
\right)^2}
<\frac{4+\varepsilon_T}{\sqrt{QT}}
$,
wo $ \varepsilon_T \to 0, \,  T\to \infty$.}

{\bf Satz 5.}
\,{\it Es sei $\alpha = (\alpha_1,\alpha_2,\alpha_3) \in
\hbox{\got S}^2 \setminus \mathbb{Q}^3$ und $\varepsilon >0$. Dann gibt es unendlich viele
rationale Vektoren $\left(\frac{A_1}{Q},\frac{A_2}{Q},\frac{A_3}{Q}\right)\in \hbox{\got S}^2\cap \mathbb{Q}^3$, so dass 
$$
\sqrt{\sum_{i=1,2,3}
\left(
\alpha_i -\frac{A_i}{Q}
\right)^2}
<
\left(
 2 \sqrt{\frac{3}{\pi}}+\varepsilon
\right)\frac{1}{Q}.
$$
}

\vskip+0.3cm
{\bf \S 3. Der Fall $n=1$. Beweis des Satzes 1 und der Bemerkung 1.}
 
Sei
$\alpha_1,\alpha_2 \in (0,1)$.
In diesem Falle folgt  aus (\ref{ret})
\begin{equation}\label{tan}
\sqrt{\sum_{i=1,2}
\left(
\alpha_i -\frac{A_i}{Q}
\right)^2} \le 
2|\arctan \beta- \arctan t|
\le \frac{2}{1+\xi^2}\, |\beta - t|
,\,\,\, t=\frac{b}{q}=\frac{b_1}{q}\in \mathbb{Q},\,\,\, \beta = \beta(\alpha )\in (0,1)
\end{equation}
(hier $\xi$ ist eine Zahl zwischen $ \beta $ und $t$).

\vskip+0.3cm
{\bf Hilfssatz 1.}\,\,{\it
Es sei $n=1$. Dann ist in (\ref{ret})
$Q=q^2+b^2$ wenn $ b\not\equiv q\pmod{2}$ und
$Q=\frac{q^2+b^2}{2}$ wenn $ b\equiv q\pmod{2}$
.

}
\vskip+0.3cm

Beweis.
Es sei  $ b\not\equiv q\pmod{2}$.
Dann $(2bq,q^2-b^2, q^2+b^2)  = 1$, sodass $Q=q^2+b^2$.

Sei  $ b\equiv q\pmod{2}$.
Dann $(2bq,q^2-b^2, q^2+b^2)  = 2$, sodass $Q=\frac{q^2+b^2}{2}$.$\Box$

\vskip+0.3cm
{\bf Hilfssatz 2. } \, {\it  
Sei $pq'-p'q\equiv 1\pmod{2}$.  Dann sind
$ p\equiv q\equiv 1\pmod{2}$, oder  $ p'\equiv q'\equiv 1\pmod{2}$, oder $ p+p'\equiv q+q'\equiv 1\pmod{2}$.

 }
\vskip+0.3cm

Beweis: Dieser Hilfssatz ist  evident.$\Box$

Beweis des Satzes 1.

{\bf Fall 1.}  $\beta$ und $\frac{1+\sqrt{5}}{2}$ seien
 \"aquivalent. 
Dann ist der Kettenbruch f\"ur  $\beta$  von der Form
 $$
\beta = [B_0; B_1,...,B_t, B_{t+1},...], \,\,\,\,\,\,\,\,\, B_\nu = 1,\,\,\nu >t. 
$$
F\"ur  den N\"aherungsbruch
$$\frac{p_\nu}{q_\nu} =[B_0;B_1,...,B_\nu]
$$
man hat
$$
\left|\beta -\frac{p_\nu}{q_\nu}\right| <\frac{1+\varepsilon}{\sqrt{5} q^2_\nu}
$$
wenn $\nu$ gro\ss \, genug ist 
[1, Kap. II].
F\"ur jedes $\nu \ge t$ liefert
Hilfssatz  2 ein $j \in \{\nu-1,\nu,  \nu+1\}$ mit $p_j \equiv q_j \pmod{2}$.
Es sei $ t_\nu = \frac{p_\nu}{q_\nu}$.
 Nun folgt f\"ur $ \frac{A_{i,\nu}}{Q_\nu } = f_i (t_\nu)$ aus (\ref{tan}) und Hilfssatz 1 
$$
Q_\nu\cdot \sqrt{\sum_{i=1,2}
\left(
\alpha_i -\frac{A_{i,\nu}}{Q_\nu}
\right)^2} 
\le \frac{q_\nu^2+p_\nu^2}{1+\xi_\nu^2}\, |\beta - t_\nu|=
  \frac{ 1+t_\nu^2}{1+\xi_\nu^2}
\cdot q_\nu\left|q_\nu\beta -{p_\nu}\right| \to \frac{1}{\sqrt{5}},
$$  
mit $\nu\to \infty$
(hier ist  $\xi_\nu$  zwischen $ \beta $ und $t_\nu$). 
Daraus folgt die Behauptung wegen
$ \sqrt{5} >\sqrt{2}$.

{\bf Fall 2.}  $\beta$ und $\frac{1+\sqrt{5}}{2}$
seien nicht  \"aquivalent. Dann
es gibt unendlich viele
 rationalen Zahlen $ t=\frac{b}{q}$ mit
$$
\left|\beta -\frac{b}{q}\right|< \frac{1}{\sqrt{8} q^2},\,\,\,\,\, (b,q) = 1
$$
(siehe [1, Kap. II]).
Hilfssatz  1 liefert $ Q\le q^2+b^2$.
 Nun folgt aus (\ref{tan}) und Hilfssatz 1 
$$
Q\cdot \sqrt{\sum_{i=1,2}
\left(
\alpha_i -\frac{A_i}{Q}
\right)^2} 
\le \frac{2(q^2+b^2)}{1+\xi^2}\, |\beta - t|=
 \frac{2(1+t^2)}{1+\xi^2}\cdot q \left|q\beta -{b}\right| \le\frac{1+\varepsilon}{\sqrt{2}},
$$  
wenn  $q $ gro\ss \, genug ist. Daraus folgt Satz 1.$\Box$

Beweis der Bemerkung 1.

Es sei $\beta =[0;2,4,\overline{2}]$. 

1. F\"ur  jeden N\"aherungsbruch
$\frac{p_\nu}{q_\nu} $
gilt $p_\nu \not\equiv q_\nu\pmod{2}$. Aber $\beta$ und $\sqrt{2}$
sind  \"aquivalent und $\frac{A_{i,\nu}}{Q_\nu} = f_i \left(\frac{p_\nu}{q_\nu}\right)$.
Nach Hilfssatz 1 ist
$$ 
Q_\nu\cdot \sqrt{\sum_{i=1,2}
\left(
\alpha_i -\frac{A_{i,\nu}}{Q_\nu}
\right)^2} 
=\frac{2(q_\nu^2+p_\nu^2)}{1+t_\nu^2+o(1)}\, |\beta - t_\nu|  \to \frac{1}{\sqrt{2}},\,\,\, \nu \to \infty.
$$ 

2. Wenn $t=\frac{b}{q}, b >9$
gilt f\"ur den Median
  $\frac{b}{q}  = \frac{p_{\nu-1}+p_{\nu}}{q_{\nu-1}+q_\nu}$  und es ist 
$ q= q_{\nu-1} +q_\nu = \frac{q_{\nu-1}+q_{\nu+1}}{2}$. Also ist
$$
\left|q\beta - b\right|
=
\frac{|q_{\nu-1}\beta - p_{\nu-1}| + |q_{\nu+1}\beta - p_{\nu+1}|}{2}
=\frac{1}{4\sqrt{2}} \left( \frac{1}{q_{\nu-1}}+\frac{1}{q_{\nu+1}}\right) (1+o(1)).
$$
Bemerken wir, dass $ q_\nu = \xi (1+\sqrt{2})^\nu (1+o(1)), \nu \to \infty $ mit positiven $\xi$.   
Aslo gilt f\"urhinreichend gro\ss es
 $q$
\begin{equation}\label{wir}
\left|\beta - \frac{b}{q}\right|
= \frac{(1 +o(1))}{\sqrt{2}q^2} .
\end{equation}

3. Sei $t=\frac{b}{q}$  weder ein N\"aherungsbruch
$\frac{p_\nu}{q_\nu} $ noch ein Median  $\frac{b}{q}  = \frac{p_{\nu-1}+p_{\nu}}{q_{\nu-1}+q_\nu}$. Dann 
folgt aus dem Satz von
 Fatou  [2,3]
 $$
\left|\beta -\frac{b}{q}\right| 
 \ge \frac{1}{q^2}.
$$
 Also   gilt  (\ref{wir})  f\"ur jeden  Bruch
 $t=\frac{b}{q}$ mit     $q$ hinreichend gro\ss\ .
Und daher haben wir
$$
Q\cdot \sqrt{\sum_{i=1,2}
\left(
\alpha_i -\frac{A_{i}}{Q}
\right)^2} 
= 2Q |\arctan \beta - \arctan t|(1-\varepsilon) \ge
(q^2+b^2) \left|\arctan \beta - \arctan \frac{b}{q}\right|(1-\varepsilon) \ge\frac{1-2\varepsilon}{\sqrt{2}}.
$$
Damit ist die
 Bemerkung bewiesen.$\Box$

\vskip+0.3cm
{\bf \S 4. Die K\"orper im $\mathbb{R}^4$.}

Es sei
$$
f({\bf w}) = f(z,y,x_1,x_2) =
zy-x_1^2-x_2^2.
$$
Wir definieren die K\"orper
$$
\hbox{\got P} =
\left\{ {\bf w} = (z,y,x_1,x_2) \in \mathbb{R}^4:\,\,
|f({\bf w})|<1 \right\}
$$
 und
$$
\hbox{\got K} =
\left\{ {\bf w} = (z,y,x_1,x_2) \in \mathbb{R}^4:\,\,
|z+y|<2,\,
|z-y|< 2\sqrt{1-\left(\frac{x_1}{\gamma}\right)^2-\left(\frac{x_2}{\gamma}\right)^2}\right\} =
$$
$$
=
\left\{ {\bf w} = (z,y,x_1,x_2) \in \mathbb{R}^4:\,\,
|\xi|<1,\,
\eta^2 +\left(\frac{x_1}{\gamma}\right)^2+\left(\frac{x_2}{\gamma}\right)^2< 1 \right\} ,
$$
wo
$$
\begin{cases}
z= \xi +\eta, \cr
y = \xi - \eta.
\end{cases}
$$
Beachte, dass $\hbox{\got K}$ konvex ist und ${\rm vol} \, \hbox{\got K} = \frac{16\pi}{3}\cdot \gamma^2 > 16$.

\vskip+0.3cm
{\bf Hilfssatz 3. } $\hbox{\got K} \subset \hbox{\got P}$.
\vskip+0.3cm

Beweis.
Es sei $ {\bf w}\in \hbox{\got K}.$ Wir haben
$
x_1^2+x_2^2 -1 < zy = \frac{(z+y)^2 -(z-y)^2}{4} < 1,
$
sodass
$ -1 < zy - x_1^2-x_2^2 < 1-x_1^2-x_2^2 \le 1$.$\Box$


Wir betrachten  die beiden
Setzen wir zwei Matrizen
$$
G_t =
\left(
\begin{array}{cccc}
t &0&0&0\cr
0&t^{-1}&0&0\cr
0&0&1&0\cr
0&0&0&1
\end{array}
\right)
,
\,\,\,\,
R_\beta
=
\left(
\begin{array}{cccc}
1&0&0&0\cr
\beta_1^2+\beta^2_2&1&-2\beta_1&-2\beta_1\cr
-\beta_1&0&1&0\cr
-\beta_2&0&0&1
\end{array}
\right),\,\,\,\,
{\rm det } \, G_t = {\rm det}\, R_\beta = 1.
$$
Es ist klar dass 
$$
f(G_t{\bf w} ) = f(R_\beta {\bf w}) = f({\bf w})
$$ 
f\"ur alle $t \in \mathbb{R}_+, \alpha \in \mathbb{R}^2, {\bf w} \in \mathbb{R}^4$.

Sei $\hbox{\got K}_\beta^t = R_\beta^{-1} G_t \hbox{\got K}$.  Aus Hilfsatz 3  folgt $\hbox{\got K}_\beta^t \subset \hbox{\got S}$,
f\"ur alle $t \in \mathbb{R}_+, \beta \in \mathbb{R}^2$. 

\vskip+0.3cm
{\bf \S 5. Beweis der S\"atze 2,3.}

Wegen ${\rm vol} \,  \hbox{\got K}_\beta^t  > 16$,
existiert nach dem Gitterpunktsatz von Minkowski
ein
Gitterpunkt 
$ {\bf g} = (q, a, b_1, b_2) \in\hbox{\got K}_\beta^t\cap \mathbb{Z}^4$, $ {\bf g} \neq{\bf 0}$.
Sei $$
L = L({\bf g})= q(\beta_1^2+\beta_2^2) +A -2b_1\beta_1 - 2b_2\beta_2,
\,\,\,
\Delta = \Delta ({\bf g} )= \sum_{i=1}^2 \left( q\beta_i -{b_i}\right)^2.
 $$
Es ist
$
R_\beta {\bf g} \in G_t \hbox{\got K}\subset \hbox{\got P}
,$
 denn
$$
|f( q, L, q\beta_1-b_1, q\beta_2 - b_2) |<1
\,\,\,\,\text{und}\,\,\,\,
|qA -b_1^2-b_2^2|<1.
$$
Aber $qA -b_1^2-b_2^2$  ist ganzzahlig, sodass
\begin{equation}\label{au}
qA = b_1^2+b_2^2,
\end{equation}
und (iii) folgt.
Es ist klar dass $q \neq 0$. Sei $q\ge 1$.
 Dann liefetr  (\ref{au}) 
 $$ L = 
\frac{q(\beta_1^2+\beta_2^2) +\frac{b_1^2+b_2^2}{q}-2b_1\beta_1 - 2b_2\beta_2}{q} =
\frac{\Delta}{q}
.$$
Wegen $ G^{-1}_t R_\beta {\bf g} =(qt^{-1}, tq^{-1}\Delta, b_1-q\beta_1, b_2 - q\beta_2)\in \hbox{\got K}$, haben wir (ii$^*$) und
$
qt^{-1} +\frac{\Delta}{qt^{-1}}
<2$.
Es folgt $\max \left(qt^{-1} , \frac{\Delta}{qt^{-1}}\right) <2$, sodass
$q < 2t$ und $ \Delta < \frac{2q}{t}$. Setzt man  $ t=T/2$,  so  haben wir  Satz 2 damit  bewiesen.

Aber  $(\beta_1, \beta_2) \not\in \mathbb{Q}^2$, so gibt es unendlich viele Vektoren    $\left(\frac{b_1}{q}, \frac{b_2}{q}\right)$  mit (ii$^*$),
und  wir haben Satz 3 bewiesen.

 \vskip+0.3cm
{\bf \S 6. \"Uber die Beweise die S\"atze 4,5.}

F\"ur $\beta =\beta (\alpha ) \in \mathbb{R}^2 $ w\'ahlen wir den Vektor
$\left(\frac{b_1}{q}, \frac{b_2}{q}\right)$  aus Satz 3 mit $\gamma \in \left(\sqrt{\frac{3}{\pi}}, \sqrt{\frac{3}{\pi}}+\frac{\varepsilon}{2}\right)$.
Wenn $ b_1^2 + b_2^2 \equiv 0 \pmod{q}$ haben wir
$
q\,|\,(2b_1q, 2b_2q, q^2-b_1^2 - b_2^2, q^2+b_1^2+b_2^2).
$
F\"ur den Nenner $Q$ aus (\ref{ret}) haben wir somit $Q\le q+\frac{b_1^2+b_2^2}{q}$ und
$$
Q\cdot
\sqrt{\sum_{i=1,2,3}  
\left(
\alpha_i -\frac{A_i}{Q}
\right)^2
}
\le q
\left(
1+\frac{b_1^2+b_2^2}{q^2}
\right)
\cdot
\frac{2}{1+r^2}
\cdot
\frac{\gamma}{q}
< 2 \sqrt{\frac{3}{\pi}}+\varepsilon
$$
(hier $r$ ist eine reelle Zahl zwischen $ \sqrt{\beta_1^2+\beta_2^2} $ und $\frac{b_1^2+b_2^2}{q^2}$).
Satz 5 folgt daraus. Der Beweis des Satzes 4
mit Hilfe von Satz 2 verl\"auft analog.
\vskip+0.3cm

Der Autor dankt  L. Summerer  für  die Hilfe in der Übersetzung.

.

\vskip+0.3cm
{\bf Literaturverzeichnis}

\noindent [1]\,
J. W. S. Cassels,\,\,
An introduction to Diophantine approximation, Cambridge University Press, 1957.

\noindent [2]\,
J. H. Grace,\,\, The classification of rational approximations,
Proc. London Math. Soc. 17 (1918), 27 - 258.

\noindent [3]\,
P. Fatou,\,\, Sur l'approximation incommencurables et les
s\'eries trigonom\'etriques,  C.R. Acad. Sci. Paris, 139 (1904), 1019 - 1021.

\noindent [4]\,
L. Fukshansky, \,\,
On similarity classes of well-rounded sublattices of $\mathbb{Z}^2$ , Journal of Number Theory 129 (2009), 2530  - 2556.

\noindent [5]\,
E. Hlawka,
\,\,
Approximation von Irrationalzahlen und pithagor\"aische Tripel,
Bonner Mathematische Schriften, 121 (1980), 1 - 32. 

\noindent [6]\,
D. Kleinbock, K. Merrill,\,\,
Rational approximation on Spheres,
preprint available at  arXiv:1301.0989v4 [math.NT] 25 May 2013.

\end{document}